\magnification=1200

\hsize=11.25cm    
\vsize=18cm       
\parindent=12pt   \parskip=5pt     

\hoffset=.5cm   
\voffset=.8cm   

\pretolerance=500 \tolerance=1000  \brokenpenalty=5000

\catcode`\@=11

\font\eightrm=cmr8         \font\eighti=cmmi8
\font\eightsy=cmsy8        \font\eightbf=cmbx8
\font\eighttt=cmtt8        \font\eightit=cmti8
\font\eightsl=cmsl8        \font\sixrm=cmr6
\font\sixi=cmmi6           \font\sixsy=cmsy6
\font\sixbf=cmbx6

\font\tengoth=eufm10 
\font\eightgoth=eufm8  
\font\sevengoth=eufm7      
\font\sixgoth=eufm6        \font\fivegoth=eufm5

\skewchar\eighti='177 \skewchar\sixi='177
\skewchar\eightsy='60 \skewchar\sixsy='60

\newfam\gothfam           \newfam\bboardfam

\def\tenpoint{
  \textfont0=\tenrm \scriptfont0=\sevenrm \scriptscriptfont0=\fiverm
  \def\rm{\fam\z@\tenrm}
  \textfont1=\teni  \scriptfont1=\seveni  \scriptscriptfont1=\fivei
  \def\oldstyle{\fam\@ne\teni}\let\old=\oldstyle
  \textfont2=\tensy \scriptfont2=\sevensy \scriptscriptfont2=\fivesy
  \textfont\gothfam=\tengoth \scriptfont\gothfam=\sevengoth
  \scriptscriptfont\gothfam=\fivegoth
  \def\goth{\fam\gothfam\tengoth}
  
  \textfont\itfam=\tenit
  \def\it{\fam\itfam\tenit}
  \textfont\slfam=\tensl
  \def\sl{\fam\slfam\tensl}
  \textfont\bffam=\tenbf \scriptfont\bffam=\sevenbf
  \scriptscriptfont\bffam=\fivebf
  \def\bf{\fam\bffam\tenbf}
  \textfont\ttfam=\tentt
  \def\tt{\fam\ttfam\tentt}
  \abovedisplayskip=12pt plus 3pt minus 9pt
  \belowdisplayskip=\abovedisplayskip
  \abovedisplayshortskip=0pt plus 3pt
  \belowdisplayshortskip=4pt plus 3pt 
  \smallskipamount=3pt plus 1pt minus 1pt
  \medskipamount=6pt plus 2pt minus 2pt
  \bigskipamount=12pt plus 4pt minus 4pt
  \normalbaselineskip=12pt
  \setbox\strutbox=\hbox{\vrule height8.5pt depth3.5pt width0pt}
  \let\bigf@nt=\tenrm       \let\smallf@nt=\sevenrm
  \normalbaselines\rm}

\def\eightpoint{
  \textfont0=\eightrm \scriptfont0=\sixrm \scriptscriptfont0=\fiverm
  \def\rm{\fam\z@\eightrm}
  \textfont1=\eighti  \scriptfont1=\sixi  \scriptscriptfont1=\fivei
  \def\oldstyle{\fam\@ne\eighti}\let\old=\oldstyle
  \textfont2=\eightsy \scriptfont2=\sixsy \scriptscriptfont2=\fivesy
  \textfont\gothfam=\eightgoth \scriptfont\gothfam=\sixgoth
  \scriptscriptfont\gothfam=\fivegoth
  \def\goth{\fam\gothfam\eightgoth}
  
  \textfont\itfam=\eightit
  \def\it{\fam\itfam\eightit}
  \textfont\slfam=\eightsl
  \def\sl{\fam\slfam\eightsl}
  \textfont\bffam=\eightbf \scriptfont\bffam=\sixbf
  \scriptscriptfont\bffam=\fivebf
  \def\bf{\fam\bffam\eightbf}
  \textfont\ttfam=\eighttt
  \def\tt{\fam\ttfam\eighttt}
  \abovedisplayskip=9pt plus 3pt minus 9pt
  \belowdisplayskip=\abovedisplayskip
  \abovedisplayshortskip=0pt plus 3pt
  \belowdisplayshortskip=3pt plus 3pt 
  \smallskipamount=2pt plus 1pt minus 1pt
  \medskipamount=4pt plus 2pt minus 1pt
  \bigskipamount=9pt plus 3pt minus 3pt
  \normalbaselineskip=9pt
  \setbox\strutbox=\hbox{\vrule height7pt depth2pt width0pt}
  \let\bigf@nt=\eightrm     \let\smallf@nt=\sixrm
  \normalbaselines\rm}

\tenpoint

\def\pc#1{\bigf@nt#1\smallf@nt}         \def\pd#1 {{\pc#1} }

\catcode`\;=\active
\def;{\relax\ifhmode\ifdim\lastskip>\z@\unskip\fi
\kern\fontdimen2  -1.2 \fontdimen3 \string;}

\catcode`\:=\active
\def:{\relax\ifhmode\ifdim\lastskip>\z@\unskip\fi\penalty\@M\ \fi\string:}

\catcode`\!=\active
\def!{\relax\ifhmode\ifdim\lastskip>\z@
\unskip\fi\kern\fontdimen2  -1.1 \fontdimen3 \string!}

\catcode`\?=\active
\def?{\relax\ifhmode\ifdim\lastskip>\z@
\unskip\fi\kern\fontdimen2  -1.1 \fontdimen3 \string?}

\frenchspacing

\def\raggedbottom{\topskip 10pt plus 36pt\r@ggedbottomtrue}

\def\pointir{\unskip . --- \ignorespaces}

\def\Medbreak{\vskip-\lastskip\medbreak}

\long\def\th#1 #2\enonce#3\endth{
   \Medbreak\noindent
   {\pc#1} {#2\unskip}\pointir{\it #3}\smallskip}

\def\decale#1{\smallbreak\hskip 28pt\llap{#1}\kern 5pt}
\def\decaledecale#1{\smallbreak\hskip 34pt\llap{#1}\kern 5pt}
\def\puce{\smallbreak\hskip 6pt{$\scriptstyle\bullet$}\kern 5pt}

\def\eqalign#1{\null\,\vcenter{\openup\jot\m@th\ialign{
\strut\hfil$\displaystyle{##}$&$\displaystyle{{}##}$\hfil
&&\quad\strut\hfil$\displaystyle{##}$&$\displaystyle{{}##}$\hfil
\crcr#1\crcr}}\,}

\catcode`\@=12

\showboxbreadth=-1  \showboxdepth=-1

\newcount\numerodesection \numerodesection=1
\def\section#1{\bigbreak
 {\bf\number\numerodesection.\ \ #1}\nobreak\medskip
 \advance\numerodesection by1}

\mathcode`A="7041 \mathcode`B="7042 \mathcode`C="7043 \mathcode`D="7044
\mathcode`E="7045 \mathcode`F="7046 \mathcode`G="7047 \mathcode`H="7048
\mathcode`I="7049 \mathcode`J="704A \mathcode`K="704B \mathcode`L="704C
\mathcode`M="704D \mathcode`N="704E \mathcode`O="704F \mathcode`P="7050
\mathcode`Q="7051 \mathcode`R="7052 \mathcode`S="7053 \mathcode`T="7054
\mathcode`U="7055 \mathcode`V="7056 \mathcode`W="7057 \mathcode`X="7058
\mathcode`Y="7059 \mathcode`Z="705A


\def\diagram#1{\def\normalbaselines{\baselineskip=0pt\lineskip=5pt}
\matrix{#1}}

\def\vfl#1#2#3{\llap{$\textstyle #1$}
\left\downarrow\vbox to#3{}\right.\rlap{$\textstyle #2$}}

\def\hfl#1#2#3{\smash{\mathop{\hbox to#3{\rightarrowfill}}\limits
^{\textstyle#1}_{\textstyle#2}}}

\def\ogoth{{\goth o}}

\def\pgoth{{\goth p}}

\def\Q{{\bf Q}}

\def\Z{{\bf Z}}

\def\F{{\bf F}}

\def\Aut{\mathop{\rm Aut}\nolimits}
\def\Hom{\mathop{\rm Hom}\nolimits}
\def\Id{\mathop{\rm Id}\nolimits}

\def\Card{\mathop{\rm Card}\nolimits}

\def\Gal{\mathop{\rm Gal}\nolimits}

\def\droite#1{\,\hfl{#1}{}{8mm}\,}

\def\to{\rightarrow}

\def\normressym(#1,#2)_#3{\displaystyle\left({#1,#2\over#3}\right)}

\def\mod{\mathop{\rm mod.}\nolimits}
\def\pmod#1{\;(\mod#1)}

\newcount\refno 
\long\def\ref#1:#2<#3>{                                        
\global\advance\refno by1\par\noindent                              
\llap{[{\bf\number\refno}]\ }{#1} \pointir{\it #2} #3\goodbreak }

\def\citer#1(#2){[{\bf\number#1}\if#2\empty\relax\else,\ {#2}\fi]}

\def\pt{\hbox{\bf.}}

\def\boxit#1{\vbox{\hrule\hbox{\vrule\kern1pt
       \vbox{\kern1pt#1\kern1pt}\kern1pt\vrule}\hrule}}
\def\cqfd{\hfill\boxit{\phantom{\i}}}

\newbox\bibbox
\setbox\bibbox\vbox{\bigbreak
\centerline{{\pc BIBLIOGRAPHY}}

\ref{\pc ALBERT} (A):
On $p$-adic fields and rational division algebras,
<Ann.\ of Math.\ (2) {\bf 41} (1940), 674--693.> 
\newcount\albert \global\albert=\refno

\ref{\pc CONRAD} (K):
Groups of order $p^3$, 
<31 August 2013, \par\vskip-.3\baselineskip\noindent
{\tt www.math.uconn.edu/\string~kconrad/blurbs/grouptheory/groupsp3.pdf}.>
\newcount\kconrad \global\kconrad=\refno

\ref{\pc DALAWAT} (C):
Further remarks on local discriminants, 
<J.\ Ramanujan Math.\ Soc.\ {\bf 25} (2010) 4, 391--417. 
Cf.~arXiv\string:0909.2541.>  
\newcount\further \global\further=\refno

\ref{\pc DALAWAT} (C):
Serre's ``\thinspace formule de masse\thinspace'' in prime degree,
<Monats\-hefte Math.\ {\bf 166} (2012), 73--92.
Cf.~arXiv\string:1004.2016.>
\newcount\monatshefte \global\monatshefte=\refno

\ref{\pc EILENBERG} (S):
Topological methods in abstract algebra. Cohomology theory of groups,
<Bull.\ Amer.\ Math.\ Soc.\ {\bf 55} (1949), 3--37.>
\newcount\eilenberg \global\eilenberg=\refno

\ref{\pc FUJISAKI} (G):
A remark on quaternion extensions of the rational $p$-adic field,
<Proc.\ Japan Acad.\ Ser.\ A Math.\ Sci.\ {\bf 66} (1990) 8, 257--259.>
\newcount\fujisaki \global\fujisaki=\refno 

\ref{\pc GREVE} (C):
Galoisgruppen von Eisensteinpolynomen {\"u}ber $p$-adischen
K{\"o}rpern, <Dissertation, Universit{\"a}t Paderborn, Oktober 2010.>
\newcount\greve \global\greve=\refno

\ref{\pc GROSS} (B) \& {\pc REEDER} (M):
Arithmetic invariants of discrete Langlands parameters,
<Duke Math.\ J.\ {\bf 154} (2010) 3, 431--508. (\$25.00)>
\newcount\grossreeder \global\grossreeder=\refno 

\ref{\pc FEIT} (W):
On $p$-regular extensions of local fields,
<Proc.\ Amer.\ Math.\ Soc.\ {\bf 10} (1959), 592--595.>
\newcount\feit \global\feit=\refno

\ref{\pc HASSE} (H):
Zahlentheorie,
<Akademie-Verlag, Berlin, 1969, 611~pp.>
\newcount\hasse \global\hasse=\refno

\ref{\pc IWASAWA} (K):
On Galois groups of local fields,
<Trans.\ Amer.\ Math.\ Soc.\ {\bf 80} (1955), 448--469.>
\newcount\iwasawa \global\iwasawa=\refno

\ref{\pc REEDER} (M):
Notes on Group Theory,
<5 October 2011,
\par\vskip-.3\baselineskip\noindent
{\tt https\string:/$\!$/www2.bc.edu/\string~reederma/Groups.pdf}.> 
\newcount\reeder \global\reeder=\refno

\ref{\pc RIOU-\pc PERRIN} (B):
Plongement d'une extension di{\'e}drale dans une extension
di{\'e}drale ou quaternionienne, 
<Th{\`e}se de troisi{\`e}me cycle (1979), Publ.\ math.\ d'Orsay, n$^\circ$~79-04.>
\newcount\perrinriou \global\perrinriou=\refno

\ref{\pc SERRE} (J-P):
Une ``\thinspace formule de masse$\,$" pour les extensions totalement
ramifi{\'e}es de degr{\'e} donn{\'e} d'un corps local, 
<Comptes Rendus {\bf 286} (1978), 1031--1036.>
\newcount\serremass \global\serremass=\refno

\ref{\pc WATERHOUSE}  (W):
The normal closures of certain Kummer extensions,
<Canad.\ Math.\ Bull.\ {\bf 37} (1994) 1, 133--139.>
\newcount\waterhouse \global\waterhouse=\refno

} 

\centerline{\bf Tame ramification and group cohomology}  
\bigskip
\centerline{Chandan Singh Dalawat} 
\centerline{Harish-Chandra Research Institute}
\centerline{Chhatnag Road, Jhunsi, Allahabad 211019, India} 
\centerline{dalawat@gmail.com}

\bigskip
\centerline{Jung-Jo Lee\footnote{*}{JJL was supported by NRF grant
No.\ 2012-005700, Republic of Korea.}}  
\centerline{Department of Mathematics, Seoul National University}
\centerline{Shillim-dong, Gwanak-gu, Seoul 151-742, Korea} 
\centerline{jungjolee@gmail.com}

\bigskip\medskip

{\bf Abstract}.  We give an intrinsic parametrisation of the set of
tamely ramified extensions of a local field with finite residue field
and bring to the fore the role played by group cohomology. We show
that two natural definitions of the cohomology class of a tamely
ramified finite galoisian extension coincide, and can be recovered
from the parameter.  We also give an elementary proof of Serre's mass
formula in the tame case and in the simplest wild case, and we
classify tame galoisian extensions of degree the cube of a prime.
\footnote{}{Keywords: 
Tame ramification,
group cohomology,
{\it 
zahme Verzweigung,
Kohomologie von Gruppen, 
ramification mod{\'e}r{\'e}e,
cohomologie des groupes
}
}

\footnote{}{MSC (2010) : 11S15, 11S20}

\bigskip\medskip

Let $K$ be a local field with finite residue field $k$ of
characteristic~$p$ and cardinality $q$.  Let $e>0$ be an integer such
that $e\not\equiv0\pmod p$ and let $f>0$ be an arbitrary integer.
Consider the set ${\cal T}_{e,f}(K)$ of $K$-isomorphism classes of
finite (separable) extensions of $K$ of ramification index~$e$ and
residual degree~$f$.  This set was investigated by Hasse in Chapter~16
of his treatise \citer\hasse(), by Albert in \citer\albert(), by
Iwasawa in \citer\iwasawa() and by Feit in \citer\feit() (sometimes
with the restriction that $K$ be of characteristic~$0$, or that the
extensions be galoisian, or that $f\not\equiv0\pmod p$).

Our purpose here is to give a more intrinsic parametrisation of this
set, and to bring to the fore the role played by group cohomology, a
theory which had not yet been fully formalised at the time of Hasse
and Albert, although only the first few cohomology groups (which were
known under different names) are needed.

We are able to recover properties of $L\in{\cal T}_{e,f}(K)$ directly
from its parameter.  These properties include those of being
galoisian, or abelian, or cyclic over $K$.  For every $L$, the
parameter determines the galoisian closure $\tilde L$ of $L$ over $K$.
When $L|K$ is galoisian, the parameter of $L|K$ determines the
cohomology class of the extension of groups
$$
1\to\Gal(L|K_f)\to\Gal(L|K)\to\Gal(K_f|K)\to1
$$ 
corresponding to the tower $L\mid K_f\mid K$, where $K_f$ is the
maximal unramified extension of $K$ in $L$~; it also determines the
smallest extension $K_{\hat f}$ of $K_f$ such that the extension of
groups corresponding to the tower $LK_{\hat f}\mid K_{\hat f}\mid K$
be split.

We also give an easy elementary proof of Serre's mass formula
\citer\serremass() in the tame case (and in the case when the
degree is divisible by $p$ but not by $p^2$), analogous to the recent
proof \citer\monatshefte() in prime degrees~$l$ (in both the cases
$l\neq p$ and $l=p$).  We explicitly work out all galoisian extensions
of $K$ of degree $l^3$ (for every prime $l\neq p$), including the case
$l=2$ of (tamely ramified) octic dihedral or quaternionic extensions.

Let $K_f$ be the degree-$f$ unramified extension of $K$,
$w_f:K_f^\times\to\Z$ its normalised valuation, $k_f$ the residue
field of $K_f$, and $G_f=\Gal(K_f|K)$.  We shall show that ${\cal
T}_{e,f}(K)$ {\it is in canonical bijection with the set of orbits for
the action of\/} $G_f$ {\it on set of what we call\/} ramified lines
$D\subset K_f^\times/K_f^{\times e}$ {\it or equivalently on the set
of sections of\/} $\bar w_f:K_f^\times/K_f^{\times e}\to\Z/e\Z$~;
ramified lines are precisely images of sections of $\bar w_f$.

We begin by recalling some basic facts about cohomology of groups
in \S{\bf 1} and apply them to the cohomology of finite fields
in \S{\bf 2}, where we verify an important compatibility between two
different definitions of the cohomology class of an extension of a
cyclic group by a cyclic group.  We then recall in \S{\bf 3} some
basic properties of the Kummer pairing such as its equivariance.  The
fundamental notion of {\it ramified lines\/} is introduced in \S{\bf
4}.  In \S{\bf 5} we parametrise the set ${\cal T}_{e,1}(K)$ and give
a proof in the spirit of \citer\monatshefte() of Serre's mass formula
in degree~$e$ (and also in degree~$ep$ when combined with the results
of \citer\monatshefte()).  We then provide in \S{\bf 6} an analogue in
degree~$e$ (prime to~$p$) of the orthogonality relation in prime
degree \citer\monatshefte().  In \S{\bf 7}, we give the
parametrisation of ${\cal T}_{e,f}(K)$ and show how the various
invariants of an $L\in{\cal T}_{e,f}(K)$ can be recovered from its
parameter.  Finally, we work out a number of instructive examples
in \S{\bf 8}.

\bigbreak
{\bf 1 Cohomology of groups}  
\medskip

Most readers can skip this \S, except perhaps ({\it 1.8\/}) where we
compute the number of $G$-orbits in a $G$-module $C$ (when both groups
ar cyclic) --- this is the key to Roquette's determination of the
cardinality of ${\cal T}_{e,f}(K)$ (7.1.4).  For an account of group
cohomology by one of its creators, see \citer\eilenberg().

\medskip
{\it 1.1 The group\/} $H^2(G,C)_\theta$
\medskip

Let $G$ be a group and $C$ a $G$-module, both written
multiplicatively, and $\theta:G\to\Aut(C)$ the action of $G$ on $C$.
An extension of $G$ by $C$ is a short exact sequence $1\to
C\to\Gamma\to G\to1$ such that the resulting conjugation action of $G$
on $C$ is equal to the given action $\theta$.  Two extensions
$\Gamma,\Gamma'$ of $G$ by $C$ are isomorphic if there is an
isomorphism of groups $\Gamma\to\Gamma'$ inducing $\Id_C$ on the
common subgroup $C$ and $\Id_G$ on the common quotient $G$.
Isomorphism classes of extensions of $G$ by $C$ are classified by the
group $H^2(G,C)_\theta$.  The class $[\Gamma]\in H^2(G,C)_\theta$
vanishes if and only if the extension $\Gamma$ is split in the sense
that the projection $\Gamma\to G$ admits a section, which happens
precisely when $\Gamma$ is isomorphic to the twisted product
$C\times_\theta G$, the product set $C\times G$ with the law of
composition $ (c,g)(d,h) =(c\theta(g)(d),gh) $.

\medbreak
{\it 1.2 The group\/} $H^1(G,C)_\theta$
\medskip

The group $H^1(G,C)_\theta$ is the set of sections of the projection
$C\times_\theta G\to G$ up to $C$-conjugacy~; it can be
identified with the set of supplements of $C$ in $\Gamma$ (subgroups
$D\subset C\times_\theta G$ such that $C\cap D=1$, $CD=\Gamma$) up to
$\Gamma$-conjugacy (or $C$-conjugacy, which comes to the same).  If
the action $\theta$ is trivial, then $H^1(G,C)_1=\Hom(G,C)$.

\medbreak
{\it 1.3 The restriction map in general}
\medskip

Let $G$ be a group and $C$ a $G$-module.  Let $\varphi:G'\to G$ be a
morphism of groups~; it allows us to view the $G$-module $C$ as a
$G'$-module via the action $\theta\circ\varphi$.  Let $C'$ be a
$G'$-module (with action $\theta'$), and let $\psi:C\to C'$ be a
morphism of $G'$-modules.  For $i=1,2$, the pair $(\varphi,\psi)$
induces a map $H^i(G,C)_\theta\to H^i(G',C')_{\theta'}$ on cohomology
called the {\it restriction map}.

For $i=1$, it sends the class in $H^1(G,C)_\theta$ of a section
$g\mapsto(\sigma(g),g)$ of the projection $C\times_\theta G\to G$ to
the class in $H^1(G',C')_{\theta'}$ of the section
$g'\mapsto(\psi(\sigma(\varphi(g'))),g')$ of the projection
$C'\times_{\theta'} G'\to G'$.

For $i=2$, the restriction map $H^2(G,C)_\theta\to
H^2(G',C')_{\theta'}$ coming from the pair $(\varphi,\psi)$ will be
defined in two steps.  In the first step, $C'=C$ and $\psi=\Id_C$, and
in the second step, $G'=G$ and $\varphi=\Id_{G'}$.

When $C'=C$ and $\psi=\Id_C$, the map $H^2(G,C)_\theta\to
H^2(G',C)_{\theta\circ\varphi}$ coming from the pair $(\varphi,\Id_C)$
sends the class of an extension $\Gamma$ of $G$ by $C$ to the class of
the extension $\Gamma_\varphi$ of $G'$ by $C$ consisting of those
$(\gamma,g')\in\Gamma\times G'$ such that $\bar\gamma=\varphi(g')$ in
$G$.

When $G'=G$ and $\varphi=\Id_{G'}$, the map $H^2(G',C)_{\theta}\to
H^2(G',C')_{\theta'}$ coming from the pair $(\Id_{G'},\psi)$ sends the
class of an extension $\Gamma'$ of $G'$ by $C$ to the class of the
extension ${}_\psi\Gamma'=(C'\times\Gamma')/\psi'(C)$ of $G'$ by $C'$,
where $\psi'(c)=(\psi(c),c^{-1})$.

The restriction map $H^2(G,C)_\theta\to H^2(G',C')_{\theta'}$ coming
from a general pair $(\varphi,\psi)$ is defined by first applying
$(\varphi,\Id_C)$ and then applying $(\Id_{G'},\psi)$ to get the
extension ${}_\psi(\Gamma_\varphi)$ of $G'$ by $C'$.
In the special case when $\varphi:G'\to G$ is surjective and $C'=C$,
the restriction map is called the {\it inflation map}~; it will be of
particular relevance in what follows.

\medbreak
{\it 1.4 The case of cyclic groups}
\medskip

Recall how the groups $H^1(G,C)_\theta$ and $H^2(G,C)_\theta$ can be
computed when $G$ is {\it cyclic\/} of order $n>0$.  Let $\sigma$ be a
generator of $G$, and define the elements $\sigma-1$ and
$N_\sigma=1+\sigma+\cdots+\sigma^{n-1}$ in the group ring $\Z[G]$
(over which $C$ is a left module via $\theta$).  We have
$N_\sigma.(\sigma-1)=0$ and $(\sigma-1).N_\sigma=0$, and therefore we
get a complex
$$
C\droite{(\ )^{\sigma-1}}
C\droite{(\ )^{N_\sigma}}
C\droite{(\ )^{\sigma-1}}
C.
\leqno{(1.4.1)}
$$
The cohomology groups of (1.4.1) are canonically isomorphic to
$H^1(G,C)_\theta$ and $H^2(G,C)_\theta$ respectively.  If $\theta$ is
trivial, then $H^1(G,C)_1={}_nC$ and $H^2(G,C)_1=C/C^n$.

Let $G'$ be another {\it cyclic\/} group, $\varphi:G'\to G$ a {\it
surjective\/} morphism of groups, and $\sigma'$ a generator of $G'$
such that $\varphi(\sigma')=\sigma$.  Let $C'$ be a $G'$-module and
$\psi:C\to C'$ a morphism of $G'$-modules.  Then the restriction map
$H^i(G,C)_\theta\to H^i(G',C')_{\theta'}$ is simply given by
restriction to subgroups and passage to the quotient from the map
$\psi:C\to C'$.

\medbreak
{\it 1.5 The case of cyclic modules}
\medskip

Specialise further to the case when $C$ is also cyclic, of some order
$m>0$, and let $a\in\Z$ be such that $\theta(\sigma)=\bar a$ in
$(\Z/m\Z)^\times$ (so that $a^n\equiv1\pmod m$).  The orders of the
cyclic groups $H^1(G,C)_a$ and $H^2(G,C)_a$ can then be computed in
terms of $a$, $m$ and $n$ because for every $r\in\Z$, the order of the
kernel ${}_rC$ (resp.~the image $C^r$) of the endomorphism $(\ )^r$ of
$C$ is $\gcd(m,r)$ (resp.~$m/\gcd(m,r)$).  Taking $r=a-1$ and
$r=1+a+\cdots+a^{n-1}$ respectively gives the result.


To get a presentation of the extension $\Gamma$ of $G$ by $C$
corresponding to a given class in $H^2(G,C)_a$, choose generators
$\tau\in C$, $\sigma\in G$ and identify the class of $\Gamma$ with the
class of an element $s\in\Z$ such that $(a-1)s\equiv0\pmod m$ (modulo
those which are ${}\equiv (1+a+\cdots+a^{n-1})t$ for some $t\in\Z$)~;
for a suitable lift $\tilde\sigma\in\Gamma$ of $\sigma$, we then have
$$
\Gamma=\langle\;\tau,\tilde\sigma\;\mid\;
\tau^m=1,\;
\tilde\sigma^n=\tau^s,\;
\tilde\sigma\tau\tilde\sigma^{-1}=\tau^a\;
\rangle.\leqno{(1.5.1)}
$$
For a direct derivation of this presentation, see for
example \citer\reeder(9.4). 

\smallbreak
{\it Example\/} 1.5.2\ \ Take $n=2$ and $m=4$.  The possibilities for
$a\pmod 4$ are $1$ and ${-1}$.  When $a=1$, we have
$H^2(G,C)_1=C/C^2$, and the two extensions $\Gamma$ (1.5.1) are the
direct product $C\times G$ and the one in which the group $\Gamma$ is
cyclic.  When $a={-1}$, we have $H^2(G,C)_{-1}={}_2C$, the split
extension is the twisted product $C\times_{-1} G$ ({\it 1.1\/}) and
called the dihedral group ${\goth D}_{4,2}$, while the other is called
the quaternionic group ${\goth Q}_{4,2}$.

\medbreak
{\it 1.6 Commutativity and cyclicity}
\medskip

Let us determine the order of $\tilde\sigma$ in $\Gamma$ (1.5.1), and
the conditions for $\Gamma$ to be commutative or cyclic.  

\smallbreak
{\it Remark\/} 1.6.1\ \ Although $s\in\Z$ is not unique in (1.5.1),
$r=\gcd(m,s)$ is uniquely determined~; $m/r$ is the order of $\tau^s$
in the group $\Gamma$.  We claim that {\it the order of the element\/}
$\tilde\sigma\in\Gamma$ (1.5.1) {\it is\/} $mn/r$.  Indeed, the order
of $\tilde\sigma$ is a multiple $dn$ of the order $n$ of its image
$\sigma\in G$~; we have to show that $d=m/r$.  Now, from the relation
$\tilde\sigma^n=\tau^s$, it follows that
$\tilde\sigma^{dn}=\tau^{ds}=1$, so $d$ is a multiple of the order
$m/r$ of $\tau^s$.  But conversely, it follows from
$\tilde\sigma^{mn/r}=\tau^{ms/r}=1$ that $dn$ divides $mn/r$ and
therefore $d$ divides $m/r$.  Hence $d=m/r$, and the order of
$\tilde\sigma$ is $mn/r$.

\smallbreak
{\it Remark\/} 1.6.2\ \  Note that the group $\Gamma$ (1.5.1) is
commutative if and only if $\tau$ and $\tilde\sigma$ commute, which
happens precisely when $a\equiv1\pmod m$, in view of the relations
$\tau^m=1$, $\tilde\sigma\tau\tilde\sigma^{-1}=\tau^a$.

\smallbreak
{\it Remark\/} 1.6.3\ \ Suppose that $a\equiv1\pmod m$.  In this case,
the extension (1.5.1) of $G$ by $C$ splits if and only if
$s\equiv0\pmod{\gcd(m,n)}$.  Indeed, this congruence is equivalent to
the existence of a $t\in\Z$ such that $nt\equiv s\pmod m$, which is
equivalent to $s\equiv(1+a+\cdots+a^{n-1})t\pmod m$ in view of
$a\equiv1\pmod m$.

\medbreak

For a prime $l$ and an integer $x\neq0$, denote by $v_l(x)$ the
exponent of $l$ in the prime decomposition of $x$.  The following
proposition has been extracted from \citer\albert(Theorem~13) and the
proof has been simplified.

\th PROPOSITION 1.6.4
\enonce
Suppose that\/ $a\equiv1\pmod m$.  The (commutative) group\/ $\Gamma$
$(1.6.2)$ is cyclic if and only if\/ $s$ is prime to\/ $\gcd(m,n)$.
\endth

{\it Proof}.  Suppose first that $s$ is prime to $\gcd(m,n)$~; we
have to find an element of order $mn$ in $\Gamma$.  The idea is to
find an element $\gamma_l\in\Gamma$ of order $l^{v_l(mn)}$ for every
prime~$l$ in each of the three (exhaustive) cases $v_l(m)v_l(n)>0$,
$v_l(m)=0$, and $v_l(n)=0$.

If $v_l(m)v_l(n)>0$, then $l$ divides $\gcd(m,n)$ and is prime to~$s$,
so $\gcd(m,s)$ is prime to~$l$ and $m/\gcd(m,s)$ is divisible by
$l^{v_l(m)}$.  Consequently, $mn/\gcd(m,s)$ is divisible by
$l^{v_l(mn)}$, and hence there is an element $\gamma_l\in\Gamma$ of
order $l^{v_l(mn)}$, in view of the fact that the order of
$\tilde\sigma\in\Gamma$ is $mn/\gcd(m,s)$ by (1.6.1).  Even if
$v_l(m)=0$ (so that $v_l(mn)=v_l(n)$), the subgroup (of order a
multiple of $n$) generated by $\tilde\sigma$ has an element $\gamma_l$
of order $l^{v_l(mn)}$.  Finally, if $v_l(n)=0$, then the subgroup (of
order $m$) generated by $\tau$ has an element $\gamma_l$ of order
$l^{v_l(m)}=l^{v_l(mn)}$.  These $\gamma_l$ are trivial for almost
all~$l$ (because $v_l(mn)=0$ for almost all~$l$), so their product
over all $l$ exists, is independent of the sequence of the factors
because $\Gamma$ is commuative by (1.6.2), and has order $mn$.

Conversely, suppose that the group $\Gamma$ is cyclic, so that
$\overline\Gamma_l=\Gamma/\Gamma^{l^{v_l(mn)}}$ (=$\Gamma\otimes\Z_l$)
is also cyclic and has order $l^{v_l(mn)}$ for every prime~$l$.
Suppose (if possible) that there is a prime~$l$ dividing all three
numbers $m,s,n$~; we shall get a contradiction by showing that
$\overline\Gamma_l$ would then have order $<l^{v_l(mn)}$.  This
follows from the fact that it is generated by the pair
$\bar\tau,\bar{\tilde\sigma}\in\overline\Gamma_l$ (images of $\tau$
and $\tilde\sigma$ respectively) each of which has order
$<l^{v_l(mn)}$, because $v_l(m)<v_l(mn)$ and
$v_l(mn/\gcd(m,s))<v_l(mn)$ by hypothesis (recall that the order of
$\tilde\sigma$ is $mn/\gcd(m,s)$ by (1.6.1)).  \cqfd

\medbreak
{\it 1.7 The inflation map in the bicyclic case}
\medskip

Let $G'$ be another cyclic group, of order $cn$ for some $c>0$, let
$\sigma'$ be a generator of $G'$, and let $\varphi:G'\to G$ be the
surjection such that $\varphi(\sigma')=\sigma$.  Regard $C$ as a
$G'$-module via $\sigma'\mapsto\sigma\mapsto(\ )^a$.  As before, the
group $H^2(G',C)_a$ can be identified with the kernel ${}_{a-1}C$ of
$(\ )^{a-1}:C\to C$ modulo the image of $(\
)^{1+a+\cdots+a^{cn-1}}:C\to C$.  Notice that
$1+a+\cdots+a^{cn-1}\equiv(1+a+\cdots+a^{n-1})c\pmod m$, 
Hence there is a commutative diagram
$$
\def\\{\mskip-2\thickmuskip}
\def\droite#1{\\\hfl{#1}{}{8mm}\\}
\diagram{
{}_{a-1}C&\droite{(\ )^c}&{}_{a-1}C\cr
\vfl{}{}{8mm}&&\vfl{}{}{8mm}\cr
H^2(G,C)_a&\droite{}&H^2(G',C)_a\cr
}\leqno{(1.7.1)}
$$
in which the vertical arrows are the passage to the quotient.  We
claim that the lower horizontal arrow --- induced by $(\ )^c$ --- is
the same as the restriction map ({\it 1.3\/}) coming from the pair
$(\varphi,\Id_C)$.

\th PROPOSITION 1.7.2
\enonce
The map\/ $H^2(G,C)_a\to H^2(G',C)_a$ in the above diagram is the
inflation map corresponding to the quotient\/ $\varphi:G'\to G$.
\endth

{\it Proof}.  Let a class in $H^2(G,C)_a$ be represented by an
extension $\Gamma$ of $G$ by $C$ having the presentation (1.5.1).  The
inflated extension $\Gamma'$ of $G'$ by $C$ consists of
$(\alpha,\beta)\in\Gamma\times G'$ such that
$\bar\alpha=\varphi(\beta)$ in $G$ ({\it 1.3}).  As a lift
$\tilde\sigma'\in\Gamma'$ of the generator $\sigma'\in G'$, we choose
$\tilde\sigma'=(\tilde\sigma,\sigma')$.  We then have
$\tilde\sigma'{}^{cn}=(\tilde\sigma^{cn},\sigma'{}^{cn})=(\tau^{cs},1)=\tau^{cs}$
and we are done, because $\Gamma'$ admits the desired presentation
$
\Gamma'=
\langle\,
\tau,\tilde\sigma'\mid
\tau^m=1,\,
\tilde\sigma'{}^{cn}=\tau^{cs},\,
\tilde\sigma'\tau\tilde\sigma'{}^{-1}=\tau^a\,
\rangle
$.
\cqfd

\medbreak
{\it 1.8 The number of orbits}
\medskip

The following lemma captures one of the basic ingredients in
Roquette's computation \citer\hasse(Chapter~16) of the number of
tamely ramified extensions of given ramification index and residual
degree (7.1.4).

\th LEMMA 1.8.1
\enonce
Let\/ $C$ be a cyclic group of order\/ $m>0$, let\/ $a>0$ be prime
to\/ $m$, and make\/ $\Z$ act on\/ $C$ by\/ $1\mapsto(\ )^a$.  The
number of orbits for this action is\/ $\sum_{t|m}\phi(t)/\chi_a(t)$,
where\/ $\chi_a(t)$ denotes the order of\/ $\bar a$ in the group\/
$(\Z/t\Z)^\times$ of order\/ $\phi(t)$.
\endth

{\it Proof}.  If $x,y\in C$ are in the same orbit, then they have the
same order in $C$.  The possible orders are the divisors of $m$~; for
each divisor $t$ of $m$, there are $\phi(t)$ elements of order
$t$. Also, the orbit of an $x\in C$ of order $t$ has $\chi_a(t)$
elements.  Indeed, if the orbit consists of the $r$
elements \hbox{$x$, $x^a$, $\ldots$, $x^{a^{r-1}}$}, then $r$ is the
smallest integer $>0$ such that $x^{a^r}=x$, or equivalently $r$ is
the smallest integer $>0$ such that $a^r\equiv1\pmod t$, so
$r=\chi_a(t)$.  \cqfd

\bigbreak
{\bf 2 Cohomology of finite fields}
\medskip

We now apply the results of \S{\bf 1} to some galoisian modules
arising from finite fields.

Let $p$ be a prime number, $k$ a finite extension of $\F_p$ with $q$
elements, $k_f$ the \hbox{degree-$f$} extension of $k$ (for every
$f>0$), and $G_f=\Gal(k_f|k)$.  Let $e>0$ be an integer such that
$q^f\equiv1\pmod e$.  We are interested in the groups
$H^2(G_f,k_f^\times\!/k_f^{\times e})_q$ and
$H^2(G_f,{}_ek_f^\times)_q$, where ${}_ek_f^\times$ is the group of
$e$-th roots of~$1$ in $k_f^\times$.  Every $\xi\in k_f^\times$ such
that $\xi^{q-1}\in k_f^{\times e}$ gives rise to a class in each of
these two $H^2$~; we prove their compatibility ({\it 2.2\/}).  For a
given class in $H^2(G_f,k_f^\times\!/k_f^{\times e})_q$, we also
determine (2.3.4) the smallest multiple $\hat f$ of $f$ such that the
inflated class vanishes in $H^2(G_{\hat f},k_{\hat f}^\times\!/k_{\hat
f}^{\times e})_q$.

\medbreak
{\it 2.1 The classes in\/} $H^2(G_f,k_f^\times\!/k_f^{\times e})_q$
{\it and\/} $H^2(G_f,{}_ek_f^\times)_q$
\medskip

Let $\xi\in k_f^\times$.  If the image $\bar\xi\in
k_f^\times\!/k_f^{\times e}$ is such that $\bar\xi^{q-1}=\bar1$, then
it has a class $[\bar\xi]\in H^2(G_f,k_f^\times\!/k_f^{\times e})_q$.
But there is also a way to attach a class in
$H^2(G_f,{}_ek_f^\times)_q$ to such $\xi$ which was inspired
by \citer\grossreeder(6.1)

Write $\xi^{q-1}=\alpha^e$ for some $\alpha\in k_f^\times$, and put
$\zeta=N_f(\alpha)$, where $N_f:k_f^\times\to k^\times$ is the norm
map.  We then have
$$
\zeta^e=N_f(\alpha^e)=N_f(\xi^{q-1})=1, 
$$
so $\zeta\in{}_e k_f^\times$.  At the same time $\zeta\in k^\times$
(being the $N_f$ of something in $k_f^\times$), so $\zeta^{q-1}=1$.
In other words, $\zeta$ is in the kernel of $(\
)^{q-1}:{}_ek_f^\times\to{}_ek_f^\times$, and so has a class
$[\zeta]\in H^2(G_f,{}_ek_f^\times)_q$.  If now we replace $\alpha$ by
$\varepsilon\alpha$ for some $\varepsilon\in{}_ek_f^\times$, then
$\zeta$ gets replaced by $N_f(\varepsilon)\zeta$.  As
$N_f(\varepsilon)=\varepsilon^{1+q+\cdots+q^{f-1}}$, the class
$[\zeta]\in H^2(G,{}_ek_f^\times)_q$ is uniquely determined by $\xi$
and does not depend on the choice of $\alpha$.

\medbreak
{\it 2.2 The compatibility of the two classes}
\medskip

Recall that $q^f\equiv1\pmod e$.  The two groups
$k_f^\times\!/k_f^{\times e}$, ${}_ek_f^{\times}$ are cyclic of the
same order $e$ and they are canonically isomorphic as $G_f$-modules
by\/ $\bar\xi\mapsto\xi^{(q^f-1)/e}$.  Therefore we get a canonical
isomorphism
$$
H^2(G, k_f^\times\!/k_f^{\times e})_q\to H^2(G,{}_ek_f^\times)_q.
$$

\th PROPOSITION 2.2.1
\enonce
Under this isomorphism, the class\/ $[\bar\xi]$ of any\/ $\xi\in
  k_f^\times$ such that\/ $\xi^{q-1}\in k_f^{\times e}$ gets mapped to
  the class\/ $[\zeta]$ of\/ $\zeta=N_f(\alpha)$ for any\/ $\alpha\in
  k_f^\times$ such that\/ $\xi^{q-1}=\alpha^e$.
\endth

{\it Proof}.  Put $S=1+q+\cdots+q^{f-1}$.  Notice first that the
condition $\xi^{q-1}\in k_f^{\times e}$ is equivalent to
$\xi^{(q^f-1)/e}\in k^\times$, because $k_f^{\times e}$
(resp.~$k^\times=k_f^{\times S}$) is the subgroup of order $(q^f-1)/e$
(resp. $q-1$) of the cyclic group $k_f^\times$ of order $q^f-1$.
Indeed, if $\omega$ is a generator of $k_f^\times$ and if
$\xi=\omega^x$, then the condition $\xi^{q-1}\in k_f^{\times e}$ is
equivalent to $x(q-1)\equiv0\pmod{e}$, and the condition
$\xi^{(q^f-1)/e}\in k^\times$ is equivalent to
$x(q^f-1)/e\equiv0\pmod{S}$.  But these two congruences are equivalent
(and are clearly satisfied when $\xi\in k^\times$~; they might
sometimes be satisfied even by some $\xi\notin k^\times$).

Now let $\xi\in k_f^\times$ be such that $\xi^{q-1}=\alpha^{e}$ for
some $\alpha\in k_f^\times$, or equivalently, as we've seen,
$\xi^{(q^f-1)/e}=\beta^{S}$ for some $\beta\in k_f^\times$.  We have to
show that $N_f(\alpha)=\alpha^{S}$ and $\beta^{S}$, which are both in
the kernel of the endomorphism $(\ )^{q-1}$ of ${}_ek_f^\times$,
define the same class in $H^2(G,{}_ek_f^\times)$ or equivalently that
$(\beta\alpha^{-1})^S=\eta^S$ for some $\eta\in{}_ek_f^\times$.

Choose a generator $\omega$ of $k_f^\times$ and write $\xi=\omega^x$,
$\alpha=\omega^a$, $\beta=\omega^b$ with $x,a,b\in\Z$, so that
$$
(q-1)x=ae+(q^f-1)c,\quad
\left({q^f-1\over e}\right)x=bS+(q^f-1)d
$$
for some $c,d\in\Z$.  We then have $(b-a)S=\left({q^f-1\over
e}\right)(cS-de)$, so if we take $\eta=\omega^{\left({q^f-1\over
e}\right)c}$, then $\eta\in{}_ek_f^\times$ and
$\eta^S=(\beta\alpha^{-1})^S$, hence $\alpha^S$ has the same class as
$\beta^S$ in $H^2(G_f,{}_ek_f^\times)$, which was to be proved.  \cqfd

\medbreak
{\it 2.3 The inflation map}
\medskip

Let $f'>0$ be a multiple of $f$.  By our notational convention,
$k_{f'}$ is the degree-$f'$ extension of $k$ and
$G_{f'}=\Gal(k_{f'}|k)$.  The inclusion $k_f^\times\to k_{f'}^\times$
induces a map on the quotients $k_f^\times\!/k_f^{\times e}\to
k_{f'}^\times\!/k_{f'}^{\times e}$.  The reader may wish to compare
the following lemma with \citer\greve(Satz~3.6).

\th LEMMA 2.3.1
\enonce
For a given\/ $\xi\in k_f^\times$, the smallest multiple\/ $f'$ of\/
$f$ such that\/ $\bar\xi^{q-1}=\bar1$ in\/
$k_{f'}^\times\!/k_{f'}^{\times e}$ is\/ $f'=df$, where\/ $d$ is the
order of\/ $\bar\xi^{q-1}$ in\/ $k_{f}^\times\!/k_{f}^{\times e}$.
\endth

{\it Proof}.  Clearly, $f'$ being a multiple of $f$, the relation
$\bar\xi^{q-1}=\bar1$ holds in\/ $k_{f'}^\times\!/k_{f'}^{\times e}$
if and only if $\xi^{q-1}\in k_{f'}^{\times e}$.  The result follows
from the fact that the degree of the extension $k_f(\root
e\of{\xi^{q-1}})$ over $k_f$ equals $d$.  \cqfd

\medbreak

Next, for every divisor $c$ of $e$, we have $k_{cf}=k_f\!\left(\root
c\of{k_f^\times}\right)$ and the natural map
$\iota:k_f^\times\!/k_f^{\times e}\to k_{cf}^\times\!/k_{cf}^{\times
e}$ is ``\thinspace raising to the exponent~$c$\thinspace'' in the
sense that if we choose a generator $\omega_c\in k_{cf}^\times$ and
put $\omega={\omega_c}^{q^{cf}-1\over q^f-1}$ (which is a generator of
$k_f^\times$), then $\iota(\bar\omega)=\bar\omega_c^c$ in
$k_{cf}^\times\!/k_{cf}^{\times e}$.  Indeed, since $q^f\equiv1\pmod
e$, we have
$$
{q^{cf}-1\over q^f-1}
=q^{(c-1)f}+\cdots+q^f+1
\equiv c\pmod e.
$$
Now, the map $\iota:k_f^\times\!/k_f^{\times e}\to
k_{cf}^\times\!/k_{cf}^{\times e}$ is $G_{cf}$-equivariant and hence
induces the inflation map
$$
H^2(G_f,k_f^\times\!/k_f^{\times e})_q\droite{}{}
 H^2(G_{cf},k_{cf}^\times\!/k_{cf}^{\times e})_q.
\leqno{(2.3.3)}
$$

\th LEMMA 2.3.4
\enonce
For a given\/ $\xi\in k_f^\times$ such that $\bar\xi^{q-1}=1$ in
$k_f^\times\!/k_f^{\times e}$, the smallest multiple\/ $\hat f$ of\/ $f$
such that\/ $[\bar\xi]=0$ in\/
$H^2(G_{\hat f},k_{\hat f}^\times\!/k_{\hat f}^{\times e})_q$ is\/ $\hat f=\hat cf$,
where\/ $\hat c$ is the order of\/ $[\bar\xi]$ in\/
$H^2(G_f,k_{f}^\times\!/k_{f}^{\times e})_q$.
\endth

{\it Proof}.  We have seen that for every divisor $c$ of $e$,
$\iota:k_f^\times\!/k_f^{\times e}\to k_{cf}^\times\!/k_{cf}^{\times
e}$ is ``\thinspace raising to the exponent~$c$\thinspace'', which is
compatible with the inflation map (2.3.3) by (1.7.2).  \cqfd

\bigbreak
{\bf 3 Kummerian extensions}
\medskip

We need to recall some basic facts about abelian extensions of
exponent dividing~$d$ of a field $F$ which contains a primitive $d$-th
root of~$1$ and which is galoisian of finite degree over some other
field $F'$.

{\it 3.1 Background}
\medskip

  Essentially as a consequence of the Hilbert-Noether vanishing
theorem for a certain $H^1$ (Satz~90), the maximal abelian extension
of $F$ of exponent dividing $d$ is $M=F(\root d\of{F^\times})$, and
there is a perfect pairing
$$
\Gal(M|F)\times(F^\times\!/F^{\times d})\longrightarrow{}_dF^\times,\quad
(\sigma,\bar x)={\sigma(y)\over y}\ \ (y^d=x)
\leqno{(3.1.1)}
$$
between the profinite group $\Gal(M|F)$ and the discrete group
$F^\times\!/F^{\times d}$.  For any closed subgroup
$H\subset\Gal(M|F)$, we have $M^H=F(\root n\of D)$ where $D\subset
F^\times\!/F^{\times d}$ is the orthogonal complement of $H$ for the
above pairing.  Conversely, for every subgroup $D\subset
F^\times\!/F^{\times d}$, the orthogonal complement
$H\subset\Gal(M|F)$ is a closed subgroup and $M^H=F\!(\root n\of D)$.
Also, for every subgroup $D\subset F^\times\!/F^{\times d}$, the
pairing (3.1.1) gives an isomorphism of (profinite) groups $\Gal(F(\!\root
d\of D)|F)\to\Hom(D,{}_dF^\times)$.

\medskip
{\it 3.2 Equivariant pairings}
\medskip

Now suppose that $F$ is a galoisian extension of finite degree over
some field $F'$, of group $G=\Gal(F|F')$.  If $D\subset
F^\times\!/F^{\times d}$ is a subgroup such that $F(\!\root d\of D)$
is galoisian over~$F'$, then the group $\Gal(F(\!\root d\of D)|F)$ may
be considered as a $G$-module for the conjugation action coming from
the short exact sequence
$$
1\to \Gal(F(\!\root d\of D)|F)
\to \Gal(F(\!\root d\of D)|F')
\to G
\to 1.
\leqno{(3.2.1)}
$$

\th PROPOSITION 3.2.2
\enonce
The extension\/ $F(\root d\of D)$ is galoisian over $F'$ if and only
if the subgroup\/ $D\subset F^\times\!/F^{\times d}$ is\/
$G$-stable. If so, the isomorphism of groups\/ $\Gal(F(\!\root d\of
D)|F)\to\Hom(D,{}_dF^\times)$ is\/ $G$-equivariant. 
\endth

{\it Proof}.  Suppose first that $D$ is $G$-stable.  We have to show
that $F(\root d\of D)$, which is clearly separable over $F'$,
coincides with all its $F'$-conjugates.  The notation $F(\root d\of
D)$ stands for $F((\root d\of x)_{x\in\tilde D})$, where $\tilde
D\subset F^\times$ is the preimage of $D$.  For $\sigma\in G$, we have
$\sigma(x)=y^dx'$ for some $x'\in \tilde D$ and some $y\in F^\times$
(because $D$ is $G$-stable), and therefore $\root
d\of{\sigma(x)}=y\root d\of{x'}$ is in $F(\root d\of D)$, so this
extension is galoisian over $F'$.

Conversely, suppose that $F(\root d\of D)$ is galoisian over $F'$, and
let $\tilde\sigma$ be an extension of some $\sigma\in G$ to an
$F'$-automorphism of $F(\root d\of D)$.  For every $x\in\tilde D$, we
have $\tilde\sigma(\root d\of x)^d=\tilde\sigma(x)=\sigma(x)$, so
$\sigma(x)\in\tilde D$ (because it has the $d$-th root
$\tilde\sigma(\root d\of x)$ in $F(\root d\of D)$), and hence $D$ is
$G$-stable.  

Finally, to check that the isomorphism $\Gal(F(\!\root d\of
D)|F)\to\Hom(D,{}_dF^\times)$ (when $D$ is $G$-stable) is $G$-equivariant, it
is enough to check that the pairing $\varphi:\Gal(F(\!\root d\of D)|F)\times
\tilde D\to{}_dF^\times$ (3.1.1) is $G$-equivariant in the sense that
$\varphi(\sigma\pt\tau,\sigma\pt x)=\sigma\pt\varphi(\tau,x)$.
Indeed, for every lift $\tilde\sigma\in\Gal(F(\root d\of D)|F')$ of a
$\sigma\in G$, we have $\tilde\sigma (\root d\of x)^d=\sigma(x)$ and
$$
\varphi(\sigma\pt\tau,\sigma\pt x)
={\tilde\sigma\tau\tilde\sigma^{-1}
 (\tilde\sigma (\root d\of x))\over\tilde\sigma(\root d\of x)} 
={\tilde\sigma\tau(\root d\of x)\over\tilde\sigma(\root d\of x)} 
=\sigma\!\left({\tau(\root d\of x)\over\root d\of x}\right)
=\sigma\pt\varphi(\tau,x)
$$
for every $\tau\in\Gal(F(\!\root d\of D)|F)$ and every $x\in\tilde
D$.  \cqfd

\smallbreak
{\it Remark\/} 3.2.3\ \ When $d$ is prime, $F'$ contains a primitive
$d$-th root of~$1$, and $G$ is a cyclic $d$-group, the class in $H^2$
of the extension (3.2.1) has been computed in \citer\waterhouse().  In
the case of interest to us, $F'$ is a local field, $F$ is finite
unramified over $F'$, $d$ is prime to the residual characteristic, and
$D$ is a $G$-stable ``ramified line'' ({\it 4.1})~; we will see later
(\S{\bf 7}) how to compute the class of (3.2.1) from $D$.

\medskip
{\it 3.3 Orbits and equivalence}
\medskip
\bigbreak

\th PROPOSITION 3.3.1
\enonce
The set of cyclic extensions of\/ $F$ of degree~$d$ up to
$F'$-isomorphisms is in natural bijection with the set of orbits for
the action of\/ $G$ on the set of cyclic subgroups of\/
$F^\times\!/F^{\times d}$ of order~$d$.
\endth

{\it Proof}.  Suppose first that the order-$d$ cyclic subgroups
$D_1,D_2\subset F^\times\!/F^{\times d}$ are in the same $G$-orbit, so
that $D_2=\sigma(D_1)$ for some $\sigma\in G$, and let $L_1=F(\root
d\of{D_1})$, $L_2=F(\root d\of{D_2})$.  Let $D_1$ be
generated by the image of $x\in F^\times$, so that $D_2$ is generated
by the image of $\sigma(x)$~; we have
$$
L_1=F[T]/(T^d-x),\quad
L_2=F[T]/(T^d-\sigma(x)).
$$
Consider the (unique) $F'$-automorphism $\tilde\sigma$ of $F[T]$ such
that $\tilde\sigma(a)=\sigma(a)$ for $a\in F$ and
$\tilde\sigma(T)=T$. Composing it with the projection $F[T]\to L_1$
induces a $F'$-morphism $L_1\to L_2$ which is an $F'$-isomorphism
because $L_1$ and $L_2$ have the same degree over~$F'$.

Conversely, if $L_i=F(\root d\of{x_i})$ for some $x_i\in F^\times$
whose images in $F^\times\!/F^{\times d}$ have order~$d$, and if we
have an $F'$-isomorphism $\tilde\sigma:L_1\to L_2$, we have to show
that $D_2=\sigma(D_1)$ for some $\sigma\in G$, where $D_i\subset
F^\times\!/F^{\times d}$ is the subgroup generated by the image of
$x_i$.  Now, $\tilde\sigma(F)=F$ because $F$ is galoisian over $F'$,
and hence $\tilde\sigma|_F=\sigma$ for some $\sigma\in G$.  Also,
$\sigma(x_1)$ has a $d$-th root in $L_2$ (namely $\tilde\sigma(\root
d\of{x_1})$) and its image has order~$d$ in $F^\times\!/F^{\times d}$,
so it generates the same subgroup as the image of $x_2$.  In other
words, $D_2=\sigma(D_1)$, and we are done.  \cqfd

\bigbreak
{\bf 4 Ramified lines}  
\medskip

Let $K$ be a local field with finite residue field $k$ of
characteristic $p$ and cardinality $q$.  Denote by $\ogoth$
(resp.~$\pgoth$) the ring of integers of $K$ (resp.~the unique maximal
ideal of $\ogoth$, so that $k=\ogoth/\pgoth$).  We have the
decomposition $\ogoth^\times=U_1.k^\times$ in which $U_1=1+\pgoth$ is
a $\Z_p$-module.  As a result, for every integer $e>0$ such that
$e\not\equiv0\pmod p$, we have the exact sequence
$$
1\to k^\times\!/k^{\times e}\to 
K^\times\!/K^{\times e}
\droite{\bar w}\Z/e\Z\to0
$$
in which $\bar w$ is induced by the normalised valuation
$w:K^\times\to\Z$.

\medbreak
{\it 4.1 The definition of ramified lines}
\medskip

The set ${\cal R}_e(K)$ of {\it ramified lines\/} consists of
subgroups $D\subset K^\times\!/K^{\times e}$ such that the restriction
$D\to\Z/e\Z$ of $\bar w$ to $D$ is an isomorphism~; ramified lines are
precisely the images of sections of $\bar w$.  As the conjugation action of
$\Z/e\Z$ on $k^\times\!/k^{\times e}$ resulting from the above exact
sequence is trivial, the number of ramified lines is equal to the
order $g=\gcd(q-1,e)$ of
$$
H^1(\Z/e\Z,k^\times\!/k^{\times e})_1
=\Hom(\Z/e\Z,k^\times\!/k^{\times e})
=k^\times\!/k^{\times e}.
$$

Every uniformiser $\pi$ of $K$ gives a bijection of the set ${\cal
R}_e(K)$ of ramified lines with the group $k^\times\!/k^{\times e}$~;
to the class $\bar u\in k^\times\!/k^{\times e}$ of $u\in k^\times$
corresponds the ramified line generated by the image of $u\pi$ in
$K^\times\!/K^{\times e}$.  Notice that the map $x\mapsto x^{q-1\over
g}$ identifies the group $k^\times\!/k^{\times e}$ with the kernel
${}_ek^\times$ of $(\ )^e:k^\times\to k^\times$.  With this
identification, to $\xi\in{}_ek^\times$ corresponds the ramified line
generated by $u\pi$ for any $u\in k^\times$ such that $u^{q-1\over
e}=\xi$.

\medbreak
{\it 4.2 The galoisian action on the set of ramified lines}
\medskip

For every $f>0$, let $K_f$ be the unramified extension of $K$ of
degree~$f$, $k_f$ its residue field, and $G_f=\Gal(K_f|K)$.  The group
$G_f$ acts on the set ${\cal R}_e(K_f)$ of ramified lines in
$K_f^\times\!/K_f^{\times e}$.  Indeed, if $D$ is generated by the
image of a uniformiser $\varpi$ of $K_f$, then $\sigma(D)$ is
generated by the image of the uniformiser $\sigma(\varpi)$ and hence
$\sigma(D)$ is a ramified line.  Also, $\Card{\cal R}_e(K_f)=g_f$,
where $g_f=\gcd(q^f-1,e)$ ({\it 4.1}).

For every uniformiser $\pi$ of $K$, the bijection
$k_f^\times\!/k_f^{\times e}\to{\cal R}_e(K_f)$ ({\it 4.1\/}) is
$G_f$-equivariant.  Therefore $\Card{\cal R}_e(K_f)^{G_f}=g$, where
$g=\gcd(q-1,e)$ is the order of ${}_{q-1}(k_f^\times\!/k_f^{\times
e})$.

\th PROPOSITION 4.2.1
\enonce
The number of orbits for the\/ $G_f$-action on\/ ${\cal R}_e(K_f)$ is
$\sum_{t|g_f}{\phi(t)/ \chi_q(t)}$ ($1.8.1$).
\endth
{\it Proof}.  Using a uniformiser of $K$, this amounts to computing
the number of orbits for the action of $G_f$ on
$k_f^\times\!/k_f^{\times e}$.  As the canonical generator of $G_f$
acts on the cyclic group $k_f^\times\!/k_f^{\times e}$ of order $g_f$
by the automorphism $(\ )^q$, the result follows from (1.8.1).  \cqfd

\medbreak
{\it 4.3 The cohomology class of a
stable ramified line, first defintion}
\medskip

Suppose that $q^f\equiv1\pmod e$ (if not, replace $e$ by
$g_f=\gcd(e,q^f-1)$).  Denote the canonical generator of $G_f$ by
$\sigma$, let $\pi$ be a uniformiser of $K$ and let $D\subset
K_f^\times\!/K_f^{\times e}$ be the ramified line generated by
$\xi\pi$ for some $\xi\in k_f^\times$.  If $D$ is $G_f$-stable, which
amounts to $\sigma(D)=D$, then $(\xi^q\pi)(\xi\pi)^{-1}\in k_f^{\times
e}$ or equivalently $\bar\xi^{q-1}=\bar1$ in $k_f^\times\!/k_f^{\times
e}$.

If we replace $\pi$ by $\pi'=u\pi$ ($u\in k^\times$), then $\bar\xi$
is replaced by $\bar\xi'=\bar\xi\bar u$.  But the norm map
$N_f:k_f^\times\to k^\times$ is surjective, so
$u=a^{1+q+\cdots+q^{f-1}}$ for some $a\in k_f^\times$, and hence
$[\xi]=[\xi']$ in $H^2(G_f,k_f^\times\!/k_f^{\times n})_q$.  Thus, the
map ${\cal R}_e(K_f)^{G_f}\to H^2(G_f, k_f^\times\!/k_f^{\times e})_q$
does not depend on the choice of $\pi$.  This defines the class $D$ in
$H^2(G_f,k_f^\times\!/k_f^{\times n})_q$.

\medbreak
{\it 4.4 The cohomology class of a stable ramified line, second
defintion}
\medskip

We assigns a class in
$H^2(G_f,{}_ek_f^\times)_q$ to $D\in{\cal R}_e(K_f)^{G_f}$ following
\citer\grossreeder(6.1).  If $D$ is
generated by $\xi\pi$, then $\xi^{q-1}=\alpha^e$ for some $\alpha\in
k_f^\times$~; put $\zeta=N_f(\alpha)$.  We have seen ({\it 2.1\/})
that $\zeta$ defines a class in $H^2(G_f,{}_ek_f^\times)_q$ which does
not depend on the choice of $\alpha$.  Moreover, if we replace $\pi$
by $\pi'=u\pi$ ($u\in k^\times$), then $\xi$ gets replaced by
$\xi'=\xi u^{-1}$, and then
$\xi'^{q-1}=\xi^{q-1}(u^{-1})^{q-1}=\alpha^e$, so we may use the same
$\alpha$ for $\pi'$ as for $\pi$.  In other words, the class $[\zeta]$
depends only on $D$.  We thus get a similar map ${\cal
R}_e(K_f)^{G_f}\to H^2(G_f,{}_ek_f)_q$.

\medbreak
{\it 4.5 The compatibility of the two definitions}
\medskip

Recall that we have an isomorphism $H^2(G_f, k_f^\times\!/k_f^{\times
e})_q\to H^2(G_f,{}_ek_f)_q$ ({\it 2.2})~; let us show that it is
compatible with the two maps from ${\cal R}_e(K_f)^{G_f}$.

\th PROPOSITION 4.5.1
\enonce
When $q^f\equiv1\pmod e$, the two definitions of the cohomology class
of\/ $D\in{\cal R}_e(K_f)^{G_f}$ are compatible with the above
isomorphism. 
\endth

{\it Proof}.  This follows from the preceding constructions and
(2.2.1). \cqfd

\medbreak
{\it 4.6 The restriction map}

Let $f'>0$ be a multiple of $f$.  If $D\in{\cal R}_e(K_f)$ is a
ramified line, generated by the image of some uniformiser $\varpi$ of
$K_f$, then the image of $\varpi$ in $K_{f'}^\times\!/K_{f'}^{\times
e}$ generates a ramified line, defining the map ${\cal
R}_e(K_f)\to{\cal R}_e(K_{f'})$. It sends $G_f$-stable ramified lines
to $G_{f'}$-stable ones, so we get the following diagram in which the
lower horizontal arrow is the restriction map (2.3.3)
$$
\def\\{\mskip-2\thickmuskip}
\def\droite#1{\\\hfl{#1}{}{8mm}\\}
\diagram{
{\cal R}_e(K_f)^{G_f}&\droite{}{}&{\cal R}_e(K_{f'})^{G_{f'}}\cr
\vfl{}{}{5mm}&&\vfl{}{}{5mm}\cr
H^2(G_f,k_f^\times\!/k_f^{\times e})_q&\droite{}{}&
 H^2(G_{f'},k_{f'}^\times\!/k_{f'}^{\times e})_q.\cr
}
\leqno{(4.6.1)}
$$

\th PROPOSITION 4.6.2
\enonce
The diagram $(4.6.1)$ is commutative.
\endth

{\it Proof}.  This follows from (1.7.2) upon choosing a uniformiser of
$K$. \cqfd

\bigbreak
{\bf 5 Totally tamely ramified extensions}.  
\medskip

Let $e>0$ be an integer $\not\equiv0\pmod p$.  Let us first study the
set ${\cal T}_{e,1}(K)$ of ($K$-isomorphism classes of) totally
ramified extensions of $K$ of degree~$e$.

\medbreak
{\it 5.1 The parametrisation of\/} ${\cal T}_{e,1}(K)$
\medskip

\th PROPOSITION 5.1.1
\enonce
The set\/ ${\cal T}_{e,1}(K)$ of totally ramified extensions of\/ $K$
of degree~$e$ is in canonical bijection with the set\/ ${\cal R}_e(K)$
of ramified lines in\/ $K^\times\!/K^{\times e}$. In particular, the
cardinality of\/ ${\cal T}_{e,1}(K)$ is $g=\gcd(q-1,e)$.
\endth

{\it Proof}.  For every uniformiser $\pi$ of $K$, the polynomial
$T^e-\pi$ is irreducible (Eisenstein's criterion) and the extension
$F(\root e\of\pi)$ is totally ramified of degree~$e$.  Conversely, let
$L|K$ be a totally ramified extension of degree~$e$ (so that the
residue field of $L$ is $k$), let $\pi$ be a uniformiser of $K$, and
write $\pi=u\xi\varpi^e$, where $u$ (resp.~$\varpi$) is a $1$-unit
(resp. uniformiser) in $L$, and $\xi\in k^\times$. Since the group of
$1$-units of $L$ is a $\Z_p$-module and $e\in\Z_p^\times$, there is a
(unique) $1$-unit $v$ of $L$ such that $u=v^e$, so the uniformiser
$\xi^{-1}\pi$ of $K$ has the $e$-th root $v\varpi$ in $L$ and
therefore $L=K(\root e\of{\xi^{-1}\pi})$.

For any two uniformisers $\pi_1,\pi_2$ of $K$, the extensions $K(\root
e\of{\pi_1})$, $K(\root e\of{\pi_2})$ are $K$-isomorphic if and only
if the unit $\pi_1/\pi_2\in\ogoth^\times$ is in $\ogoth^{\times e}$,
which happens precisely when $\pi_1$ and $\pi_2$ generate the same
ramified line in $K^\times\!/K^{\times e}$.  This completes the
proof.  \cqfd

\th PROPOSITION 5.1.4
\enonce
For\/ $L\in{\cal T}_{e,1}(K)$, the group\/ $\Aut_K(L)$ is canonically
isomorphic to\/ ${}_ek^\times$ and hence it is cyclic of order\/
$g=\gcd(q-1,e)$.
\endth

{\it Proof}.  Indeed, $L=K(\root e\of\pi)$ for some uniformiser $\pi$
of $K$, and the $K$-conjugates of $\root e\of\pi$ in $L$ are precisely
$\xi\root e\of\pi$, where $\xi$ is an $e$-th root of~$1$ in $K$.  The
map $\sigma\mapsto\sigma(\root e\of\pi)/\root e\of\pi$ is thus an
isomorphism $\Aut_K(L)\to{}_ek^\times$.

This isomorphism is independant of the choice of $\pi$.  Indeed, every
other uniformiser $\pi'$ of $K$ such that $L=K(\root e\of\pi')$ is
of the form $\pi'=\varepsilon^e\pi$ for some $\varepsilon\in
k^{\times}$ (ignoring $1$-units of $K$, which we can).  We may thus
take $\varepsilon\root e\of\pi$ for $\root e\of\pi'$, and we have
$$
{\sigma(\root e\of{\pi'})\over\root e\of{\pi'}}
={\sigma(\varepsilon\root e\of\pi)\over\varepsilon\root e\of\pi}
={\varepsilon\sigma(\root e\of\pi)\over\varepsilon\root e\of\pi}
={\sigma(\root e\of\pi)\over\root e\of\pi}
$$
for every $\sigma\in\Aut_K(L)$, which was to be proved. \cqfd

\th COROLLARY 5.1.5
\enonce
Some\/ $L\in{\cal T}_{e,1}(K)$ is galoisian over\/ $K$ if and only if
$e$ divides $q-1$.  If so, then every\/ $L\in{\cal T}_{e,1}(K)$ is
galoisian (and indeed cyclic) over\/ $K$.
\endth

{\it Proof}.  A finite separable extension $L$ of $K$ is galoisian
over $K$ if and only if $\Aut_K(L)$ has order $[L:K]$.  For an
$L\in{\cal T}_{e,1}(K)$, this happens precisely when $\gcd(q-1,e)=e$
(5.1.4), or equivalently when $e$ divides $q-1$.  \cqfd

\medbreak
{\it 5.2 Serre's mass formula in tame degrees}
\medskip

For the next corollary, we need to recall the statement of Serre's
mass formula \citer\serremass().  Let $n>0$ be any integer and denote
by ${\cal T}_{n,1}(K)$ the set of $K$-isomorphism classes of finite
(separable) totally ramified extensions of $K$ of ramification
index~$n$.  For every $L\in{\cal T}_{n,1}(K)$, put
$c_K(L)=w(\delta_{L|K})-(n-1)$, where $\delta_{L|K}$ is the
discriminant of $L|K$.  The mass formula asserts that
$$
\sum_{L\in{\cal T}_{n,1}(K)}{1\over|\Aut_K(L)|}q^{-c_K(L)}=n.
\leqno{(5.2.1)}
$$ 
where $|\Aut_K(L)|$ is the order of the group of $K$-automorphisms of
$L$. 

\th COROLLARY 5.2.2
\enonce
Serre's mass formula\/ $(5.2.1)$ holds over\/ $K$
in every tame degree\/~$e$ (prime to\/ $p$).
\endth

{\it Proof}.  Indeed, for every $L\in{\cal T}_{e,1}(K)$, we have
$c_K(L)=0$, $|\Aut_K(L)|=g$ (5.1.4) and there are $g$ such $L$
(5.1.1), where $g=\gcd(q-1,e)$. \cqfd

In fact we can do slightly better if we use the results
of \citer\monatshefte() where a new proof of Serre's mass formula in
degree~$p$ was given.  Let $\tilde K$ be a separable algebraic closure
of $K$, and let $E\subset\tilde K$ run through totally ramified
extensions of degree~$n$ over $K$, which we express by $[E]\in{\cal
T}_{n,1}(K)$.  Serre \citer\serremass() shows that $(5.2.1)$ is equivalent
to
$$
\sum_{E\subset\tilde K,\;[E]\in{\cal T}_{n,1}(K)} q^{-c_K(E)}=n.
\leqno{(5.2.3)}
$$ 

\th PROPOSITION 5.2.4
\enonce
Serre's mass formula\/ $(5.2.3)$ holds over\/ $K$ in degree\/~$n=ep$
(with\/ $e\not\equiv0\pmod p$).
\endth

{\it Proof}.  Let $E\subset\tilde K$ be a totally ramified extension
of degree $ep$ over $K$, and let $L$ be the maximal tamely ramified
extension of $K$ in $E$~; we have $[L:K]=e$.  By the formula for the
transitivity of the discriminant, we have
$$
w(\delta_{E|K})=(e-1)p+w_L(\delta_{E|L})
$$
where $w$ (resp.~$w_L$) is the normalised valuation of $K$
(resp.~$L$).  It follows that $c_K(E)=c_L(E)$.  Next, notice that
there are precisely $e$ totally ramified extensions of $K$ in $\tilde
K$ of degree $e$ over~$K$, since there are $g=\gcd(q-1,e)$ isomorphism
classes in ${\cal T}_{e,1}(K)$ (5.1.1), and each class $[L]$ is
represented by $e/g$ extensions $L\subset\tilde K$, because
$g=|\Aut_K(L)|$ (5.1.4).  Now, by decomposing the sum
$\sum_{[E:K]=ep}$ (5.2.3) as $\sum_{[L:K]=e}\sum_{[E:L]=p}$, we have
$$
\sum_{[E:K]=ep} q^{-c_K(E)}
=\sum_{[L:K]=e}\sum_{[E:L]=p}q^{-c_K(E)}
=\sum_{[L:K]=e}\sum_{[E:L]=p}q^{-c_L(E)}.
$$
But $\sum_{[E:L]=p}q^{-c_L(E)}=p$ by \citer\monatshefte(th.~35), and
hence $\sum_{[E:K]=ep} q^{-c_K(E)}=ep$, as was to be proved. \cqfd

\smallbreak
{\it Remark\/} 5.2.5\ \ The same d{\'e}vissage reduces the proof of
(5.2.3) for arbitray~$n$ to the case $n=p^r$.  Note that a proof of
(5.2.1) for $n$ prime can also be found in \citer\monatshefte().

\bigbreak
{\bf 6 The orthogonality relation}  
\medskip

Let us make some remarks about the special case $q\equiv1\pmod e$.
More precisely, suppose that the (cyclic) group ${}_eK^\times\subset
K^\times$ of $e$-th roots of~$1$ in $K$ has order~$e$, and let
$M=K(\root e\of{K^\times})$ be the maximal abelian extension of $K$ of
exponent dividing $e$.  We have the perfect pairing (3.1.1)
$$
\Gal(M|K)\times(K^\times\!/K^{\times e})\longrightarrow{}_eK^\times
$$
of free rank-$2$ $(\Z/e\Z)$-modules, defined by $(\sigma,\bar
x)=\sigma(y)/y$ for any $y\in M^\times$ such that $y^e=x$.  

\th PROPOSITION 6.0.1
\enonce
The orthogonal complement of the inertia subgroup\/ $\Gamma_0$ of\/
$\Gal(M|K)$ is the subgroup\/ $k^\times\!/k^{\times e}$ of\/
$K^\times\!/K^{\times e}$ (and conversely).
\endth

{\it Proof}.  Indeed, the fixed field $M^{\Gamma_0}$ of the inertia
subgroup is the maximal unramified extension $M_0$ of $K$ in $M$.  It
is easy to see that, $\omega$ being a generator of $k^\times$, the
extension $K(\root e\of\omega)$ of $K$ in $M$ is unramified and of
degree~$e$ over $K$.  At the same time, the ramification index of
$M|K$ is at least $e$, as it contains $K(\root e\of\pi)$ for any
uniformiser $\pi$ of $K$.  As $[M:K]=e^2$, we must have $M_0=K(\root
e\of{k^\times})$, which was to be proved. \cqfd

\th PROPOSITION 6.0.2
\enonce
For every subgroup\/ $D\subset K^\times\!/K^{\times e}$, the maximal
unramified extension\/ of\/ $K$ in\/ $K(\root e\of D)$ is\/ $K(\root
e\of{D_0})$, with $D_0=D\cap (k^\times\!/k^{\times e})$.
\endth

{\it Proof}.  Let $L=K(\root e\of D)$, and let $L_0$ be the maximal
unramified extension of $K$ in $L$~; it is clear that $K(\root
e\of{D_0})\subset L_0$.  Conversely, if $C\subset k^\times\!/k^{\times
e}$ is the subgroup such that $L_0=K(\root e\of C)$, then $C\subset D$
and hence $C\subset D_0$.  It follows that $C=D_0$. \cqfd

\smallbreak
{\it Remark\/} 6.0.3\ \ As a corollary, $L=K(\root e\of{D})$ is
totally ramified over $K$ if and only if $D\cap(k^\times\!/k^{\times
e})=\{1\}$.  The analogue of (6.0.1) in degree~$p$ can be found
in \citer\further() (\S1) if $K$ has characteristic~0 and
in \citer\further() (\S5) if $K$ has characteristic~$p$.

\bigbreak
{\bf 7 The parametrisation of\/} ${\cal T}_{e,f}(K)$   
\medskip

Recall that ${\cal T}_{e,f}(K)$ is the set of $K$-isomorphism classes
of separable extensions of $K$ of ramification index~$e$
($\not\equiv0\pmod p$) and residual degree~$f$.  We will see that it
can be identified with the set of orbits for the action of
$G_f=\Gal(K_f|K)$ on the set ${\cal R}_e(K_f)$ of ramified lines in
$K_f^\times\!/K_f^{\times e}$.

There is a canonical surjection ${\cal T}_{e,1}(K_f)\to{\cal
T}_{e,f}(K)$, and a canonical bijection ${\cal T}_{e,1}(K_f)\to{\cal
R}_e(K_f)$ (by (5.1.1), applied to $K_f$), so the question is~: When
are the extensions defined by two distinct ramified lines in
$K_f^\times\!/K_f^{\times e}$ isomorphic as extensions of $K$
(although they are not $K_f$-isomorphic)~?

\medbreak
{\it 7.1 The parametrisation of\/} ${\cal T}_{e,f}(K)$
\medskip

\th PROPOSTION 7.1.1
\enonce
The extensions\/ $L,L'$ corresponsing to two ramified lines\/
$D,D'\subset K_f^\times\!/K_f^{\times e}$ are\/ $K$-isomorphic if and
only if\/ $D'=\sigma(D)$ for some\/ $\sigma\in G_f$.
\endth

{\it Proof}.  The proof is similar to that of (3.3.1), although there
the extensions $L,L'$ were kummerian whereas here they need not even
be galoisian (over $K_f$).  Suppose first that $D'=\sigma(D)$, and let
$\varpi$ be a uniformiser of $K_f$ whose image generates $D$, so that the
image of $\sigma(\varpi)$ generates $D'$, and
$$
L=K_f[T]/(T^e-\varpi),\quad
L'=K_f[T]/(T^e-\sigma(\varpi)).
$$
Consider the (unique) $K$-automorphism $\tilde\sigma$ of $K_f[T]$ such
that $\tilde\sigma(a)=\sigma(a)$ for every $a\in K_f$ and
$\tilde\sigma(T)=T$. Composing it with the projection $K_f[T]\to L'$
induces a $K$-morphism $L\to L'$ which is a $K$-isomorphism.

Conversely, if $L=K_f(\root e\of\varpi)$ for some uniformiser $\varpi$ of
$K_f$, and if we have a $K$-isomorphism $\tilde\sigma:L\to L'$, then
its restriction to the maximal unramified extensions of $K$ in $L$ and
$L'$ is a $K$-automorphism $\sigma:K_f\to K_f$, and the uniformiser
$\sigma(\varpi)$ of $K_f$ has the $e$-th root $\tilde\sigma(\root
e\of\varpi)$ in $L'$, so $L'=K_f(\root e\of{\sigma(\varpi)})$.  In other
words, $D'=\sigma(D)$.  \cqfd

\th COROLLARY 7.1.2
\enonce
The set\/ ${\cal T}_{e,f}(K)$ is in natural bijection with the set of
orbits\/ ${\cal R}_e(K_f)/\!\!/G_f$ for the action of\/ $G_f$ on\/
${\cal R}_e(K_f)$.\cqfd
\endth

\th COROLLARY 7.1.3
\enonce
An extension\/ $L\in{\cal T}_{e,1}(K_f)$ is galoisian over\/ $K$ if
and only if the corresponding\/ $G_f$-orbit consists of a single\/
$D\in{\cal R}_e(K_f)$ and\/ $q^f\equiv1\pmod e$.
\endth

{\it Proof}.  Indeed, for $L$ to be galoisian over $K$ it must be
galoisian over $K_f$, which is equivalent to $q^f\equiv1\pmod e$
(5.1.5), and all $K$-conjugates of $L$ must coincide, which is
equivalent to $D\in{\cal R}_e(K_f)^{G_f}$ (3.2.2).  \cqfd

\smallbreak
{\it Remark\/} 7.1.4\ \ It follows from the parametrisation (7.1.2)
that the set\/ ${\cal T}_{e,f}(K)$ has
$\sum_{t|g_f}{\phi(t)/\chi_q(t)}$ elements, in the notation of
(4.2.1), where $g_f=\gcd(q^f-1,e)$.  If $q^f\equiv1\pmod e$ (in which
case $g_f=e$), precisely $g=\gcd(q-1,e)$ of these are galoisian
over\/~$K$, by (7.1.3) and ({\it 4.2}).  If $q\equiv1\pmod e$ (in
which case $g_f=g=e$), the $G_f$-action on the set ${\cal R}_e(K_f)$
is trivial, so ${\cal T}_{e,f}(K)$ contains $e$ extensions and all of
them are abelian over\/~$K$ (1.6.2).  These are the only galoisian or
abelian cases. Cf.~\citer\hasse(Chapter~16)

\smallbreak
{\it Remark\/} 7.1.5\ \ For $L\in{\cal T}_{e,f}(K)$ galoisian of group
$G=\Gal(L|K)$ and inertia subgroup $G_0$, the short exact sequence
$1\to G_0\to G\to G/G_0\to1$ splits if and only if $L=K_f(\root
e\of\pi)$ for some uniformiser $\pi$ of $K$.  Indeed, suppose first
that $L=K_f(\root e\of{\pi})$, and let $E=K(\root e\of\pi)$.  As $E$
is totally ramified of degree~$e$ over $K$, the extension $L$ of $E$
is unramified (and hence cyclic) of degree~$f$~; it can be seen that
$\Gal(L|E)$ is a supplement of $G_0$ in $G$.  Conversely, if $G_0$ has
a supplement $S$ in $G$, then the extension $L^S$ of $K$ is totally
ramified of degree~$e$ and hence of the form $K(\root e\of\pi)$ for
some uniformiser $\pi$ of $K$, and $L=K_f(\root e\of\pi)$.  \cqfd

\medbreak
{\it 7.2 The presentation of the group}
\medskip

Suppose that $q^f\equiv1\pmod e$ and let $D\in{\cal R}_e(K_f)^{G_f}$,
so that the extension $L=K_f(\root e\of D)$ is in ${\cal T}_{e,f}(K)$
and galoisian over $K$ (7.1.3).  The inertia subgroup
$\Gamma_0=\Gal(L|K_f)$ of $\Gamma=\Gal(L|K)$ is canonically isomorphic
to $\Hom(D,{}_eK_f^\times)={}_eK_f^\times$ (because $D$ is isomorphic
to $\Z/e\Z$ by $\bar w$), or more simply by (5.1.4)), and the
identification $\Gamma_0={}_eK_f^\times$ is $G_f$-equivariant (3.2.2).
We thus have an extension
$$
1\to{}_eK_f^\times\to\Gamma\to G_f\to1
\leqno{(7.2.1)}
$$
and we would like to compute its class in $H^2(G_f,{}_eK_f^\times)_q$
in terms of the parameter $D$ of $L$.  

\th PROPOSITION 7.2.2
\enonce
The class of the extension\/ $(7.2.1)$ is the same as the class\/ $[D]\in
H^2(G_f,{}_ek_f^\times)_q$ of\/ $D$ $(4.5.1)$.
\endth

{\it Proof}.  We will actually compute a presentation of the group
$\Gamma$ (as in \citer\grossreeder()) and observe that it is the
extension corresponding to the class of $D$ as defined in ({\it
4.4\/}).

Let $\pi$ be a uniformiser of $K$ and suppose that the $G_f$-stable
ramified line $D$ is generated by (the image of) $\xi\pi$ for some
$\xi\in k_f^\times$ (such that $\xi^{q-1}=\alpha^e$ for some
$\alpha\in k_f^\times$), so that $L=K_f(\root e\of{\xi\pi})$.  

Choose a generator $\tau$ of $\Gamma_0$, so that $\tau(\root
e\of{\xi\pi})=\zeta\root e\of{\xi\pi}$ for a certain (generator)
$\zeta\in{}_eK_f^\times$.  Notice that $N_f(\xi)^{q-1}=1$, so
$N_f(\alpha)^e=1$, and hence $N_f(\alpha)=\zeta^s$ for some $s\pmod
e$.  As $N_f(\alpha)\in k^\times$, we must have $(q-1)s\equiv0\pmod
e$.

Also choose a lift $\tilde\sigma\in\Gamma$ of the
canonical generator $\sigma\in G_f$.  Now,
$$
\tilde\sigma(\root e\of{\xi\pi})^e=
\tilde\sigma(\xi\pi)=
\xi^{q-1}.\xi\pi=
(\alpha\root e\of{\xi\pi})^e,
$$
so that $\tilde\sigma(\root e\of{\xi\pi})=\zeta^j\alpha\root
e\of{\xi\pi}$ for some $j\pmod e$.  Replacing $\tilde\sigma$ by
$\tau^{-j}\tilde\sigma$, we may assume that $\tilde\sigma(\root
e\of{\xi\pi})=\alpha\root e\of{\xi\pi}$.  We then have
$\tilde\sigma^2(\root e\of{\xi\pi})=\tilde\sigma(\alpha)\alpha.\root
e\of{\xi\pi}$ and so on, hence
$$
\tilde\sigma^f(\root e\of{\xi\pi})
=N_f(\alpha)\root e\of{\xi\pi}
=\zeta^s\root e\of{\xi\pi}
=\tau^s(\root e\of{\xi\pi}).
$$
It follows that $\tilde\sigma^f=\tau^s$.  Finally, 
$$
\tau^q\tilde\sigma(\root e\of{\xi\pi})
=\tau^q(\alpha\root e\of{\xi\pi})
=\zeta^q\alpha\root e\of{\xi\pi}
=\tilde\sigma(\zeta\root e\of{\xi\pi})
=\tilde\sigma\tau(\root e\of{\xi\pi}),
$$
and hence $\tilde\sigma\tau{\tilde\sigma}^{-1}=\tau^q$.  We have found
that the group $\Gamma$ (of order $ef$) is generated by
$\langle\tau,\tilde\sigma\rangle$, and the relations
$$
\tau^e=1,\ \ \ 
\tilde\sigma^f=\tau^s,\ \ \ 
\tilde\sigma\tau{\tilde\sigma}^{-1}=\tau^q
$$
hold.  But we have seen that the group (1.5.1) with this presentation
has $ef$ elements, so this is indeed a presentation for $\Gamma$.  So
the class of $D$ is the same as the class of the extension
(7.2.1).  \cqfd

\medbreak
{\it 7.3 The invariants of an orbit}
\medskip

We are now going to review a certain number of invariants of a
$G_f$-orbit in ${\cal R}_e(K_f)$ which recover the invariants of the
corresponding $L\in{\cal T}_{e,f}(K)$ such as the galoisian closure
$\tilde L$ of $L$ over $K$, or the smallest extension $K_{\hat f}$ of
$K_f$ for which the exact sequence
$1\to\Gamma_0\to\Gamma\to\Gamma/\Gamma_0\to1$ splits, where
$\Gamma=\Gal(LK_{\hat f}|K)$ and $\Gamma_0$ is the inertia subgroup of
$\Gamma$.  

(7.3.1) In general, let $L\in{\cal T}_{e,f}(K)$, and let $\tilde L$ be
the galoisian closure of $L$ over $K$.  It is clear that $\tilde L$ is
tamely ramified over $K$, so $\tilde L\in{\cal T}_{\tilde e, cf}(K)$
for some multiple $\tilde e$ of $e$ and some $c>0$.  As $\tilde e$
(and hence $e$) divides $q^{cf}-1$ (7.1.3), $c$ is a multiple of the
order $r$ of $q^f$ in $(\Z/e\Z)^\times$, and therefore
$K_{rf}\subset\tilde L$.  Replacing $L$ by $LK_{rf}$, we assume that
$q^f\equiv1\pmod e$.

Let $D\in{\cal R}_e(K_f)$ be a ramified line representing the
$G_f$-orbit corresponding to $L$ and let $\xi\in k_f^\times$ be such
that $D$ is generated by $\xi\pi$ (so that $L=K_f(\root
e\of{\xi\pi})$).  The order $d$ of $\bar\xi^{q-1}$ in
$k_f^\times\!/k_f^{\times e}$ depends only on $L$, not on the choices
of $D$ and $\pi$, and the galoisian closure of $L$ over $K$ is $\tilde
L=K_{df}(\root e\of{\xi\pi})$ (2.3.1).  In particular, $\tilde e=e$.

Indeed, if we replace $\pi$ by $\pi'=u\pi$ for some $u\in k^\times$
and $D$ by $\sigma(D)$ for some $\sigma\in G_f$, then $\bar\xi$ gets
replaced by $\sigma(\bar\xi\bar u^{-1})$.  But then $\bar\xi^{q-1}$
and  $\sigma(\bar\xi\bar u^{-1})^{q-1}$ have the same order because
$\sigma$ is an automorphism of $k_f^\times\!/k_f^{\times e}$ and
$u^{q-1}=1$.  

(7.3.2) Now suppose that $L\in {\cal T}_{e,f}(K)$ is galoisian over
$K$.  What is the smallest $\hat c>0$ such that the extension $\hat
L=LK_{\hat cf}$ splits over $K$ in the sense that $\hat L=K_{\hat
cf}(\root e\of\pi)$ for some uniformiser $\pi$ of $K$~? This is
equivalent to the extension $\Gal(\hat L|K)$ of $G_{\hat cf}$ by
${}_eK_{\hat cf}^\times$ being split.  Now, the class in
$H^2(G_f,{}_eK_f^\times)_q$ of the extension $
1\to{}_eK_f^\times\to\Gal(L|K)\to G_f\to1 $ is the same as the class
of its parameter $D\in{\cal R}_e(K_f)^{G_f}$ (4.5.1), (7.2.2), and
hence $\hat c$ is the order of this class (2.3.4).

\bigbreak
{\bf 8 Examples}
\medskip

Recall the notation in force~: $K$ is a local field with finite
residue field $k$ of characteristic~$p$ and cardinality~$q$.  For
$f>0$, $K_f$ is the unramified extension of $K$ of degree~$f$, $k_f$
is its residue field, and $G_f=\Gal(K_f|K)$.  In order to write down
extensions of $K$ explicitly, we choose a uniformiser $\pi$ of $K$ and
a compatible system of generators $\omega_f$ of the cyclic groups
$k_f^\times$.  For $e>0$ such that $e\not\equiv0\pmod p$, ${\cal
T}_{e,f}(K)$ is the set of $K$-isomorphism classes of extensions of
$K$ of ramification index $e$ and residual degree~$f$.  The choice of
$\pi$ allows us to identify ${\cal T}_{e,f}(K)$ with the set of orbits
for the action of $G_f$ on $k_f^\times\!/k_f^{\times e}$.

 We compute all ${\goth S}_3$-extensions of $K$ ($p\neq3$), all {\it
tame\/} ${\goth S}_3$-extensions of $K$ ($p=3$) and, for every prime
$l\neq p$, all galoisian extensions of $K$ of degree $l^3$ which are
not abelian over $K$.  We also analyse all extensions $L$ in ${\cal
T}_{3,2}(K)$ ($p\neq3$) (8.2) or ${\cal T}_{4,2}(K)$ ($p\neq2$) (8.6)
by determining their galoisian closure $\tilde L$ over $K$ and the
smallest $\hat f$ such that $\tilde LK_{\hat f}$ splits over $K$ in
the sense of (7.3.2).

\th PROPOSITION 8.1
\enonce
If\/ $q\equiv-1\pmod3$, then\/ $K$ has a unique ${\goth
S}_3$-extension, namely\/ $K(\root3\of1,\root3\of\pi)$.  If\/
$q\equiv1\pmod3$, then $K$ has no\/ ${\goth S}_3$-extensions, and if\/
$p=3$, then\/ $K$ has no\/ {\rm tamely ramified\/} ${\goth
S}_3$-extensions.
\endth

{\it Proof}.  Let $L$ be an ${\goth S}_3$-extension of $K$.  If
$p\neq3$, we have $(e,f)=(3,2)$ (so $L$ is tame even when $p=2$)~: for
in all other cases ${\goth S}_3$ would have to have a quotient of
order~$3$.  A similar reasoning shows that if $p=3$, then
$e\equiv0\pmod3$, making $L$ wildly ramified over $K$.

So assume that $p\neq3$.  If $q\equiv1\pmod3$, then every $L\in {\cal
T}_{3,2}(K)$ is abelian over $K$, so $K$ doesn't have any ${\goth
S}_3$-extensions.  If $q\equiv-1\pmod3$, then the only extension in
${\cal T}_{3,2}(K)$ which is galoisian is
$L=K(\root3\of1,\root3\of\pi)$, and $\Gal(L|K)={\goth S}_3$.  \cqfd

\smallbreak
{\it Remark\/} 8.2\ \ When $p=3$, ${\goth S}_3$-extensions of $K$
correspond bijectively to separable cubic extensions which are not
cyclic over~$K$~; they are classified in \citer\monatshefte().  (More
generally, for any $p$, all separable extensions of degree~$p$ over
$K$ are parametrised, and the ones which are cyclic have been
characterised).  Suppose now that $p\neq3$.  If $q\equiv1\pmod3$, then
${\cal T}_{3,2}(K)$ consists of three extensions, all three abelian
(in fact cyclic) and split over $K$.  If\/ $q\equiv-1\pmod3$, then
${\cal T}_{3,2}(K)$ consists of two extensions, the ${\goth
S}_3$-extension $K(\root3\of1,\root3\of\pi)$ and the extension
$L=K_2(\root 3\of{\omega_2\pi})$ which is not galoisian over $K$.  The
galoisian closure of $L$ over $K$ is $\tilde L=K_6(\root3\of\pi)$
which is split over $K$.  The special case $K=\Q_2$, $\pi=2$ is
treated in \citer\greve(Beispiel~3.1).

(8.3)\ \ Let $l$ be a prime. Recall that there are exactly two groups
$\Gamma$ of order $l^3$ which is not commutative~; see for
example \citer\kconrad().  The centre $Z\subset\Gamma$ of both these
groups has order~$l$, and the quotient $\Gamma/Z$ is commutative of
exponent~$l$ (and order $l^2$).  For $l=2$, they are the dihedral
group ${\goth D}_{4,2}$ and the quaternionic group ${\goth Q}_{4,2}$
(1.5.2).  When $l\neq2$, one $\Gamma$ has exponent~$l$ (the
``Heisenberg group'' ${\goth H}_{l^3}$) and the other has
exponent~$l^2$.  The latter is the twisted product
$(\Z/l^2\Z)\times_\iota(U_1/U_2)$, where $\iota$ is the natural action
of $\Aut(\Z/l^i\Z)=(\Z/l^i\Z)^\times=\Z_l^\times\!/U_i$ (with
$U_j=1+l^j\Z_l$). We denote this group by ${\goth D}_{l^2,l}$.  (More
generally, one has the twisted product ${\goth
D}_{l^n,l^{n-r}}=(\Z/l^n\Z)\times_\iota(U_r/U_n)$ of order $l^{2n-r}$
for every $n>0$ and every $r\in[1,n]$.)

\th LEMMA 8.4
\enonce
If\/ $K$ has a galoisian extension of\/ degree\/ $l^3$ (\/$l\neq p)$
which is not abelian, then\/ $(e,f)=(l^2,l)$ and $v_l(q-1)=1$.
\endth

{\it Proof}.  If $K$ has such an extension $L$, then $K$ has an
abelian extension of degree $l^2$ and exponent $l$ (8.1), so we must
have $v_l(q-1)>0$ (7.2.4). Next, we must have $(e,f)=(l^2,l)$ because
$L|K$ is not abelian.  For the same reason, $q\not\equiv1\pmod{l^2}$,
so $v_l(q-1)=1$. \cqfd

\th
PROPOSITION 8.5
\enonce
If\/ $p\neq2$, then $K$ has a ${\goth D}_{4,2}$-extension or a\/
${\goth Q}_{4,2}$-extension if and only if\/ $q\equiv-1\pmod4$.  If so,
$K$ has a unique ${\goth D}_{4,2}$-extension and a unique\/ ${\goth
Q}_{4,2}$-extension.
\endth

{\it Proof}. If $K$ has an extension of degree~$2^3$ which is
galoisian but not abelian, then we must have $v_2(q-1)=1$ (or
equivalently $q\equiv-1\pmod4$) and $(e,f)=(4,2)$, by (8.4).

Suppose that $q\equiv-1\pmod4$.  There are three orbits for the action
of $G_2$ on $k_2^\times\!/k_2^{\times4}$, namely $\{1\}$,
$\{\bar\omega_2^2\}$, and $\{\bar\omega_2,\bar\omega_2^{-1}\}$.  So
there are two extensions in ${\cal T}_{4,2}(K)$ which are galoisian
over $K$, namely $L^{(0)}=K_2(\root4\of\pi)$ and
$L^{(2)}=K_2(\root4\of{\omega_2^2\pi})$.  Of these, $L^{(0)}$ is split
over $K$, so $\Gal(L^{(0)}|K)$ is the dihedral group~${\goth D}_{4,2}$
(1.5.2), whereas $L^{(2)}$ is not split over $K$, so $\Gal(L^{(2)}|K)$
is the quaternionic group~${\goth Q}_{4,2}$ (1.5.2).  This concludes the
proof. \cqfd

\smallbreak
{\it Remark\/} 8.6\ \ An explicit generation of the ${\goth
Q}_{4,2}$-extension when $K=\Q_p$ (and $p\equiv-1\pmod4$) can be found in
\citer\fujisaki().  
Let us analyse the set ${\cal T}_{4,2}(K)$ when $p\neq2$.  If
$q\equiv1\pmod4$, then it consists of {\it four\/} extensions, and all
four are abelian over $K$, but only two of them are split over~$K$~;
the other two (which are cyclic) split in ${\cal T}_{4,4}(K)$.  If
$q\equiv-1\pmod4$, then ${\cal T}_{4,2}(K)$ has only {\it three\/}
extensions, only two of which are galoisian over $K$, and only one of
them (the ${\goth D}_{4,2}$-extension, namely
$K(\root4\of1,\root4\of\pi$)) is split~; the other (the ${\goth
Q}_{4,2}$-extension) splits in ${\cal T}_{4,4}(K)$.  The galoisian
closure of the third $L\in{\cal T}_{4,2}(K)$ is $\tilde L=LK_4$, and
$\tilde L$ splits only in ${\cal T}_{4,8}(K)$.

\th PROPOSITION 8.7
\enonce
If\/ $l\neq2,p$, then\/ $K$ has a
galoisian extension\/ $L$ of degree\/~$l^3$ which is not abelian if
and only if\/ $v_l(q-1)=1$. If so, there are\/ $l$ such extensions\/
$L$, and the group\/ $\Gal(L|K)$ is isomorphic to\/ ${\goth
D}_{l^2,l}=(\Z/l^2\Z)\times_\iota(U_1/U_2)$ $(8.3)$ for each\/
$L$.
\endth

{\it Proof}.  If $K$ has such an extension, then $v_l(q-1)=1$ and
$(e,f)=(l^2,l)$ (8.4).  Conversely, suppose that $v_l(q-1)=1$.
Extensions in ${\cal T}_{l^2,l}(K)$ which are galoisian over $K$
correspond to the fixed points for the action of $G_l$ on
$k_l^\times\!/k_l^{\times l^2}$.  This group is cyclic of order $l^2$
because $q^l\equiv1\pmod{l^2}$.  As $\gcd(q-1,l^2)=l$, there are $l$
fixed points, so there are $l$ extensions $L\in {\cal T}_{l^2,l}(K)$
galoisian over $K$~; none of them is abelian over $K$ because $l^2$
does not divide $q-1$.  For each such $L$, the group
$\Gamma=\Gal(L|K)$ has order $l^3$ and contains the cyclic subgroup
$\Gal(L|K_l)$ of order $l^2$, so $\Gamma$ is isomorphic to ${\goth
D}_{l^2,l}$ (8.3).  The same conclusion can also be arrived at by
showing that $H^2(G_l,k_l^\times\!/k_l^{\times l^2})_q$
vanishes. \cqfd

\th COROLLARY 8.8 \enonce For\/ $l\neq2$, the Heisenberg
group\/ ${\goth H}_{l^3}$ $(8.1)$ does not occur as\/ $\Gal(L|K)$ for
any local field\/ $K$ of residual characteristic\/~$p\neq
l$. \cqfd \endth

The reader may wish to analyse the set ${\cal T}_{l^2,l}(K)$ in the
same way as we analysed ${\cal T}_{2^2,2}(K)$ in (8.6).

\smallbreak
{\it Example\/} 8.9\ \ Consider the case $q\equiv1\pmod{2^2}$.  We
have seen that every galoisian extension in ${\cal T}_{2^2,2}(K)$ is
in fact abelian.  But for some $m>2$, there might be galoisian
extensions in ${\cal T}_{2^m,2}(K)$ which are not abelian.  A
necessary and sufficient condition for that to happen is that $2^m$
divide $q^2-1$ but not $q-1$.  In view of $v_2(q+1)=1$, this condition
is equivalent to $v_2(q-1)=m-1$.

When $v_2(q-1)=m-1$, there are $2^{m-1}$ extensions $L\in {\cal
T}_{2^m,2}(K)$ which are galoisian but not abelian~; for every such
$L$, the resulting short exact sequence (7.2.1)
$$
1\to{}_{2^m}K_2^\times\to\Gal(L|K)\to G_2\to1
$$
splits because the group $H^2(G_2,k_2^\times\!/k_2^{\times 2^m})_q$
vanishes.  For some related results, see \citer\perrinriou(I.2).  

\smallbreak
{\it Remark\/} 8.10\ \ It is possible to determine all galoisian
extensions of $K$ of degree~$l^n$ (for any prime $l\neq p$ and any
$n>0$) by fixing $f=l^b$ and considering $e=l^a$ (such that $a+b=n$).
Feit \citer\feit() counts the number of $G$-extensions of $K$ when
$G$ has order prime to~$p$~; one should be able to recover his results
from the foregoing.

\bigbreak
{\bf 9 Acknowledgements} One of the authors (CSD) wishes to express
his heartfelt thanks to Rakesh Yadav, Georg-August-Universit{\"a}t,
G{\"o}ttingen, for supplying \citer\fujisaki(), and to M. Franck
Pierron, Universit{\'e} Paris-Sud, Orsay, for
supplying \citer\perrinriou().

\bigbreak
\unvbox\bibbox 

\bye